\documentclass[11pt]{article}
\usepackage{amssymb}
\usepackage{latexsym}
\usepackage{epsfig}

\newtheorem{theorem}{Theorem}
\newtheorem{lemma}{Lemma}
\newtheorem{corollary}{Corollary}
\newtheorem{definition}{Definition}
\newtheorem{remark}{Remark}

\newcommand{\la}{\lambda}
\newcommand{\lie}{\hat{\mathfrak g}}

\title{On multi-color partitions and the generalized \\Rogers-Ramanujan
identities}

\author{
Naihuan Jing
\thanks{Research supported in part by NSA grant MDA 904-97-1-0062 and
NSF grant DMS-9701755 at MSRI. }\\
{\small Department of Mathematics }\\
{\small North Carolina State University}\\
{\small Raleigh, North Carolina 27695-8205, USA}\\
{\small \tt jing@math.ncsu.edu}
\and
Kailash C. Misra
\thanks{Research supported in part by NSA grant MDA 904-00-1-0042. } \\
{\small Department of Mathematics }\\
{\small North Carolina State University}\\
{\small Raleigh, North Carolina 27695-8205, USA}\\
{\small \tt misra@math.ncsu.edu}
\and
Carla D. Savage
\thanks{Research supported in part by NSF grant
DMS9622772 and NSA grant MDA 904-00-1-0059}\\
{\small Department of Computer Science}\\
{\small North Carolina State University}\\
{\small Raleigh, North Carolina 27695-8206, USA}\\
{\small \tt savage@cayley.csc.ncsu.edu}
}

\begin{document}
\date{November 29, 2000}
\maketitle
\begin{abstract}
Basil Gordon, in the sixties, and George Andrews, in the seventies,
generalized the
Rogers-Ramanujan identities to higher moduli. These identities
arise in many areas of mathematics and mathematical physics. One of these
areas is representation theory of infinite dimensional Lie algebras,
where
various known interpretations of these identities have led to interesting
applications. Motivated by their connections  with Lie algebra
representation theory, we give a new interpretation of a sum related
to
generalized Rogers-Ramanujan identities in terms
of multi-color partitions.
\end{abstract}

\section{Introduction}
The celebrated  Rogers-Ramanujan identities and their generalizations
(see \cite{G}, \cite{A1})
have influenced current research in many areas
of mathematics
and physics (see \cite{A2,BeM}).
These identities can be expressed as:
\begin{equation}
\label{RR1}
\prod_{n\not \equiv  0, \pm 2 (\bmod
5)} (1-q^n)^{-1}\ \ =\ \ \sum_{n\geq 0}
\frac{q^{n^2}}{(1-q)(1-q^2) \ldots (1-q^n)}
\end{equation}
and
\begin{equation}
\label{RR2}
\prod_{n\not \equiv 0, \pm 1 (\bmod
5)}(1-q^n)^{-1}\ \ =\ \ \sum_{n\geq 0}
\frac{q^{n^2+n}}{(1-q)(1-q^2) \ldots (1-q^n)}.
\end{equation}

Identities (\ref{RR1}) and (\ref{RR2}) have a natural combinatorial
interpretation in terms of partitions,  which was generalized
by Gordon (\cite{G}, \cite{A1}, Theorem 7.5).
A partition of a positive integer $n$ is a finite, non-increasing sequence
of positive
integers, called parts, whose sum is $n$.
\begin{theorem}{\rm [Gordon]}
\label{gordon}
For $M=2k+1$ and $0<r \leq k$, the number of partitions of $n$ of the
form
$(\pi_1, \pi_2, \dots, \pi_l)$, where $\pi_j - \pi_{j+k-1} \geq 2$
and at most $r-1$ of the parts are 1 is equal to the number of
partitions of $n$ into parts not congruent to 0, $r$, or $-r$
modulo $M$.
\end{theorem}
Setting $r=2$ and $M=5$ in Theorem \ref{gordon} gives (\ref{RR1})
and $r=1$, $M=5$  gives (\ref{RR2}).

About twenty years ago it was observed that these
identities
play
an important role in the representation theory of affine Lie algebras
via its principal characters \cite{LM}.
In 1978 Lepowsky and Wilson \cite{LW0} gave the first explicit
realization
of the affine Lie algebra $\widehat{\mathfrak {sl}}(2)$.
This led to a new algebraic structure called the Z-algebra
\cite{LW1} which gave a formal foundation to study systematically the
connection between affine Lie algebras and combinatorial identities.
In particular, Lepowsky and Wilson \cite{LW1, LW1a, LW2} used
the Z-algebra structure to construct the integrable highest weight
representations of the affine Lie algebra
$\widehat{\mathfrak{sl}}(2)$ and gave a $Z$-algebraic proof
of the Rogers-Ramanujan identities and $Z$-algebraic interpretation of the
generalized Rogers-Ramanujan identities.

One can also use (generalized) Rogers-Ramanujan identities
to construct explicitly integrable representations of other
affine Lie algebras.
See for example \cite{BM}, \cite{M1}, \cite{M2}, \cite{Ma}, \cite{X}
for these developments. In this connection the $Z$-operators still
play an important role.
These operators act on a certain space $\Omega(V)$ called the vacuum space
associated with the representation space $V$ \cite{LW1a} and are
parameterized by the set
of roots of the associated simple Lie algebra ${\mathfrak g}$.
However,
on $V$ many of these $Z$-operators  are scalar multiples of each other. For
example, let us consider the affine  Lie algebra
$\lie=\widehat{\mathfrak{sl}}(5)$ and its integrable highest weight
representation $V(\lambda)$,  with highest weight $\lambda=
\Lambda_0+\Lambda_2$, where $\Lambda_i$ are the fundamental weights of
the Lie algebra $\lie$. The principal character of $V(\lambda)$ is

$$
\chi(V(\la))\ \ =\ \ F\prod_{n\geq 1,\\
n\not \equiv 0, \pm 3\,(\bmod 7)}(1-q^n)^{-1},
$$
where $F=\prod_{n\geq 1,\\n\not \equiv 0 \bmod 5} (1-q^n)^{-1}$. In this case
there
are
two independent families of
$Z$-operators:
$$
Z(\beta, z)=\sum_{i\in\mathbb Z}Z(\beta, i)z^{-i}, \qquad Z(\beta,
i)\in {\rm End} V(\la)
$$
for $\beta=\beta_1$ and $\beta_1+\beta_2$, where $\{\beta_1, \beta_2,
\beta_3, \beta_4\}$ are the simple roots of ${\mathfrak{sl}}(5)$
corresponding to the principal Cartan subalgebra $\mathfrak a$ (see \cite{M3}).
In
\cite{M3} $V(\Lambda_0+\Lambda_2)$ has been constructed using only
one set of $Z$-operaters $Z(\beta_1, z)$ and Gordon's
generalization of the Rogers-Ramanujan identities
with $r=3$ and $M=7$:

\begin{equation}
\label{r3M7}
\prod_{n\not \equiv 0, \pm 3 \,(\bmod 7)}(1-q^n)^{-1}\ \ =\ \ \sum_{n\geq
0}a(n)q^n,
\end{equation}
where $a(n)$ denotes the number of partitions of $n$ such that the outer
two of any three
consecutive parts differ by at least 2 and at most two parts are
1. However, from the representation theory
point of view it would be more natural to construct the representation using
both families of operators $Z(\beta_1, z)$ and $Z(\beta_1+\beta_2,
z)$.   It is expected that this would
correspond to another expansion of the left-hand side
of (\ref{r3M7}), namely
\begin{equation}\label{genRR}
\prod_{n\not \equiv 0, \pm 3 (\bmod
7)}(1-q^n)^{-1}\ \ =\ \ \sum_{n_1\geq  n_2\geq
0}\frac{q^{n_1^2+n_2^2}}{(q)_{n_1-n_2}(q)_{n_2}},
\end{equation}
where we let $(a)_n=\prod_{k=0}^{n-1}(1-aq^k)$.
The expansion (\ref{genRR}) is a special case of Andrews' and Bressoud's
analytic generalization of the Rogers-Ramanujan identities:
\begin{equation}
\prod_{n\not \equiv 0, \pm r (\bmod
2k+s)}(1-q^n)^{-1}\ \ =\ \
\sum_{n_1 \geq  \ldots \geq n_{k-1}\geq 0}
\frac{q^{n_1^2+n_2^2+ \ldots + n_{k-1}^2 + n_r + n_{r+1} + \ldots + n_{k-1}}}
{(q)_{n_1-n_2} \ldots (q)_{n_{k-2}-n_{k-1}}(q^{2-s}; q^{2-s})_{n_{k-1}}},
\label{verygenRR}
\end{equation}
where
$s=0, 1$.
Andrews
(\cite{A1}, Theorem 7.8) derived this generalization for the case
of odd modulus ($s=1$) and Bressoud \cite{B1} for the case of even
modulus ($s=0$).
In \cite{FS} Feigin and Stoyanovsky used the representation of
\cite{LP} in the homogeneous gradation to give
certain combinatorial interpretations of the multisum side
of the generalized Rogers-Ramanujan identities. Later Georgiev \cite{Ge}
and also Meurman and Primc (\cite{MP}, \cite{MP1}, \cite{MP2} and \cite{P})
related the sum sides of various generalized Rogers-Ramanujan type expressions
to multi-color partitions by attaching colors to different roots in $Z$-algebra
type constructions of the
homogeneous irreducible highest weight representations of certain affine Lie
algebras. It is clear from this work that the language of multi-color
partitions are suitable for Z-algebraic constructions and interpretations of
Rogers-Ramanujan type identities.
The Z-operator constructions in the principal gradation
motivated us to seek a new multi-color interpretation
(corresponding to the parameters $n_1,n_2, \ldots, n_{k-1}$
in the sum side of (\ref{verygenRR}))
of the generalized
Rogers-Ramanujan identities.

It is well-known that the product side of
(\ref{verygenRR})
can also be written as
\begin{equation}
\label{other_sum}
\prod_{n\not \equiv 0, \pm r (\bmod
M)}(1-q^n)^{-1} \ \ = \ \
\frac1{(q)_{\infty}}\sum_{j=-\infty}^{\infty}
(-1)^jq^{j[(M)j+M-2r]/2}
\end{equation}
a result which follows from the Jacobi Triple Product Identity.
In this paper,
we provide a combinatorial description
of the sum side of (\ref{other_sum})
in terms of multi-color partitions ($k-1$ colors
for identity (\ref{verygenRR}))
with the hope that this will
give new insights into the Z-operator constructions in the principal gradation.
Our work builds on an interpretation of the
sum side of (\ref{other_sum})
in terms
of partitions with bounds on successive ranks due to
Andrews \cite{A3} and Bressoud \cite{B3}.
Andrews and Bressoud showed that the
sum side of (\ref{other_sum}) is the generating function for
$|A_n(M,r)|$, where
$A_n(M,r)$ is the set of all partitions of $n$ whose successive
ranks lie in the interval $[-r+2,M-r-2]$ (to be discussed in more
detail in Section 2).
Our main theorem establishes a bijection between $A_n(M,r)$ and
a family of  multi-color partitions as described below.

\begin{definition}\label{D:2}
For $t \geq 1$, a {\em $t$-color partition} of $n$ is a pair
$(\alpha, c_{\alpha})$  where
$\alpha = (\alpha_1, \alpha_2, \ldots, \alpha_l)$
is a partition of $n$ and
$c_{\alpha}$ is a function which assigns to each $i \in \{1,2, \ldots, l\}$
one of the colors $\{1,2, \ldots, t\}$
so that if $\alpha_i = \alpha_{i+1}$ then
$c_{\alpha}(i) \leq c_{\alpha}(i+1)$.
We say that $c_{\alpha}(i)$ is the {\em color} of the $i$-th part of $\alpha$.
\end{definition}
For example, $(8_2,8_3,5_1,4_1,4_2,4_3,3_2,2_1)$ is a 3-color partition of 38,
where the
subscript of a part denotes its color.

Our main theorem
is stated below and
proved in Section 2 as Theorem \ref{main}.

\noindent{\bf Main Theorem} \it For integers $r$, $M$, and $k$ satisfying
$0 < r \leq M/2$ and
$k = \lfloor M/2 \rfloor$,
let $C_n(M,r)$ be the set of
$k-1$-color partitions
$(\alpha,c_{\alpha})$ of $n$
satisfying the
following three conditions.
Let
$\alpha = (\alpha_1, \alpha_2, \ldots, \alpha_l)$.

{(i) (\rm Initial Conditions)}
For $ 1 \leq i \leq l$,
\[
\alpha_i   > \left\{ \begin{array}{ll}
	|2c_{\alpha}(i)-r+1| & \mbox{{\rm if} $\alpha_i \equiv r (\bmod{2})$}\\
	|2c_{\alpha}(i)-r| & \mbox{{\rm otherwise}}
	\end{array}
	\right. .
\]

{\rm (ii) (Color Difference Conditions)}
For $1 \leq i < l$,
\[
\alpha_i - \alpha_{i+1} \geq \left\{ \begin{array}{ll}
	2 + |2(c_{\alpha}(i) - c_{\alpha}(i+1))| & \mbox{{\rm if} $\alpha_i
\equiv
		\alpha_{i+1} (\bmod{2})$}\\
	2 + |2(c_{\alpha}(i) - c_{\alpha}(i+1))-1| & \mbox{{\rm if} $\alpha_i
\not\equiv
		\alpha_{i+1} \equiv r (\bmod{2})$}\\
	2 + |2(c_{\alpha}(i) - c_{\alpha}(i+1))+1| & \mbox{{\rm if}
$\alpha_{i+1} \not\equiv
		\alpha_{i} \equiv r (\bmod{2})$}
	\end{array}
        \right. .
\]

{\rm (iii) ( Parity Condition on Last Color when $M$ is Even)}
For $1 \leq i \leq l$, if $M$ is even and $c_{\alpha}(i) = k-1$,
the last color, then
\[
\alpha_i \not\equiv r \pmod 2.
\]
Then
\[
\sum_{n=0}^{\infty}|C_n(M,r)|q^n \ \  =
\sum_{n=0}^{\infty}|A_n(M,r)|q^n \ \  =
\frac1{(q)_{\infty}}\sum_{j=-\infty}^{\infty}
(-1)^jq^{j[(M)j+M-2r]/2}
\]

\rm


Our results suggest several interesting problems.
Our motivating problem is to make use of the multi-color
partition interpretation to
construct natural realizations of integrable representations of affine Lie
algebras, via a correspondence between the $k-1$ colors
and the parameters $n_1, \ldots, n_{k-1}$ of (\ref{verygenRR}).
However, another intriguing problem
is to establish
a direct
bijection between the multi-color partitions and the partitions
counted by the sum side of
the Rogers-Ramanujan identity
(\ref{verygenRR}), possibly through $A_n(M,r)$.
We note some related work for the sum side.
In \cite{A4}, Andrews gave a combinatorial interpretation of the sum
side of (\ref{verygenRR}) in terms of Durfee dissection partitions.
Another result, due to Burge \cite{Burge1,Burge2} and formulated in terms of
lattice paths by Bressoud \cite{B4}, interprets the sum side of
(\ref{verygenRR}) as the number of lattice paths of weight $n$
starting at $(0,k-r)$ which have no peak of height $k$ or greater.
(Steps allowed in the lattice path are:  $(x,y) \rightarrow (x+1,y+1)$;
$(x,0) \rightarrow (x+1,0)$;
and, if $y>0$,
$(x,y) \rightarrow (x+1,y-1)$.
The weight of a lattice path is the sum of the $x$-coordinates
of its peaks.)

Although we have never seen a multi-color interpretation of
(\ref{verygenRR}), together with its conditions, made
explicit, multi-color interpretations of other identities of
the Rogers-Ramanujan type appear, for example, in
\cite{AgAn}, \cite{AgBr}, \cite{Ge}, (\cite{MP}, \cite{MP1}, \cite{MP2} and
\cite{P}). In Section 3 we note that ideas implicit in the
papers \cite{AgAn} and  \cite{AgBr} give rise to an alternative multi-color
interpretation of (\ref{verygenRR}).

In Section 4, we give an example of a context in which the
multi-color interpretation seems quite natural.
In \cite{FQ}, Foda and Quano derive a
finitization of a form of the generalized Rogers-Ramanujan
identities. We show that the corresponding refinement for our multi-color
partitions is simply an additional constraint on the size
of the largest part.

\section{Multi-color partitions}

In this section, we give a new combinatorial interpretation of
the product side of (\ref{verygenRR}).
Our main tool will be a combinatorial generalization
of the Rogers-Ramanujan identities, due to George Andrews,
which involves the {\em successive ranks} of a partition.

A partition
$\pi=(\pi_1, \pi_2, \dots, \pi_l)$
can be visualized by its Ferrers diagram, an array of dots, where
$\pi_i$ is the number of left justified dots in the $i$th row. The largest
square subarray of dots in this diagram is the Durfee square and the Durfee
square size, denoted by $d(\pi)$, is the length of a side. Flipping the diagram
along its main diagonal, one obtains the dual diagram, associated with its
dual partition $\pi'=(\pi_1', \pi_2', \dots, \pi_{\pi_1}')$, where
$\pi_i'$ is the number of indices $j$ with $\pi_j \geq i$.
The sequence of successive ranks of $\pi$ is the
sequence  $(\pi_1 - \pi_1', \pi_2 - \pi_2', \dots, \pi_d - \pi_d')$, where $d =
d(\pi)$.

Let $\pi$ be a partition of $n$ with successive ranks
$r_1, r_2, \ldots ,r_{d}$.
For $1 \leq i \leq d= d(\pi)$, let
\[
\alpha_i = \pi_i + \pi'_{i} -2(i-1) -1  =\pi_i + \pi'_i -2i + 1,
\]
i.e., $\alpha_i$ is the number of dots on the $i$-th ``angle'' of $\pi$.
Denote by $\alpha(\pi)$ the partition with parts:
$\alpha_1$, $\alpha_2$, $\cdots$, $\alpha_d$.
Note that $\alpha(\pi)$ is a partition of $n$ such that
\begin{equation}
\alpha_i \equiv r_i+1 \pmod 2.
\label{rankparity}
\end{equation}
(See Figure 1.)

\begin{figure}
\centering\psfig{figure=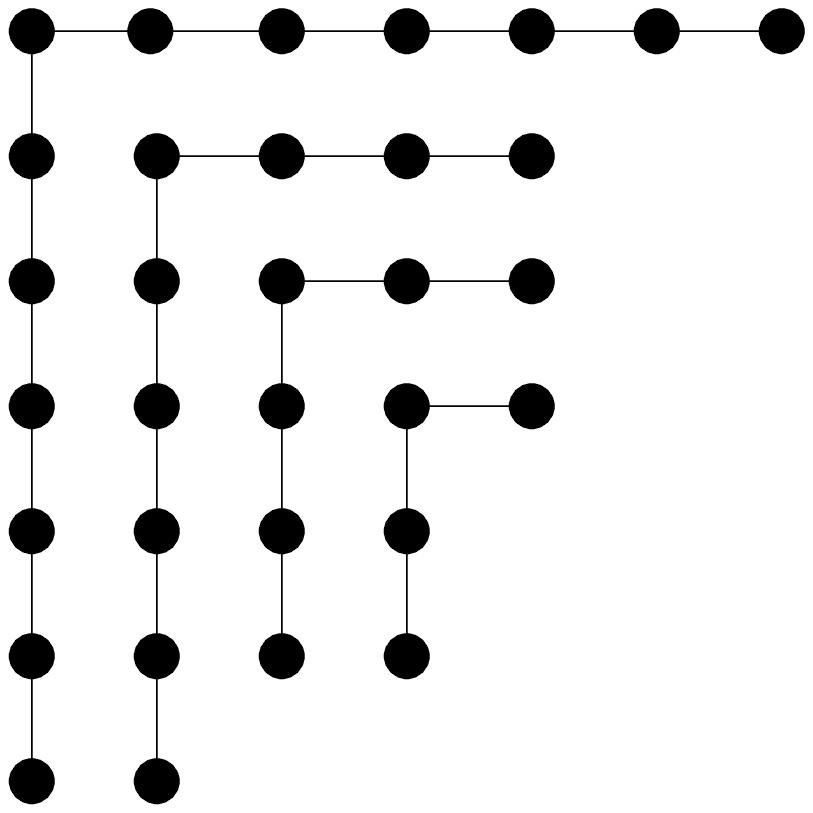,width=0.2\linewidth}
\caption{The Ferrers diagram of $\pi = (7,5,5,5,4,4,2)$ with successive ranks
$0,-2,-1,-1$, and the four angles
indicated, illustrating that
$\alpha(\pi)=(13,9,6,4)$.}
\label{angles}
\end{figure}

\begin{lemma}
For $1 \leq i \leq d(\pi)$, $\alpha_i > |r_i|$.
\label{partsize}
\end{lemma}

\noindent
{\bf Proof.}
Since $i \leq d(\pi)$, we have
$\pi_i \geq i$
and
$\pi'_i \geq i$, so
\[
\alpha_i \   =  \  \pi_i + \pi'_i -2i +1
           \ =  \  r_i + 2\pi'_i -2i+1
           \ \geq  \  r_i + 2i-2i+1
           \ =  \  r_i + 1
\]
and
\[
\alpha_i \   =  \  \pi_i + \pi'_i -2i +1
           \ =  \  -r_i + 2\pi_i -2i+1
           \ \geq  \  -r_i + 2i-2i+1
           \ =  \  -r_i + 1
\]
$\Box$

\begin{lemma}
\label{diff}
For $1\leq  i < d(\pi)$,
\[ \alpha_i - \alpha_{i+1}\  \geq\  2 + |r_i - r_{i+1}|.
\]
\end{lemma}

\noindent
{\bf Proof.}
If $r_i \geq r_{i+1}$, then since $\pi'_i \geq \pi'_{i+1}$,
\[
\pi_i - \pi_{i+1}\  \geq\
\pi_i - \pi_{i+1}
-(\pi'_i - \pi'_{i+1}) \ = \ r_i - r_{i+1}.
\]
Thus
\[ \alpha_i - \alpha_{i+1}\
 =
\pi_i - \pi_{i+1} +
\pi'_i - \pi'_{i+1} + 2
\  \geq\
r_i - r_{i+1} + 2.
\]
Similarly,
If $r_i \leq r_{i+1}$, then since $\pi_i \geq \pi_{i+1}$,
\[
\pi'_i - \pi'_{i+1}\  \geq\
\pi'_i - \pi'_{i+1}
-(\pi_i - \pi_{i+1}) \ = \ -r_i + r_{i+1}.
\]
Thus
\[ \alpha_i - \alpha_{i+1}\
 =
\pi'_i - \pi'_{i+1} +
\pi_i - \pi_{i+1} + 2
\  \geq\
-r_i + r_{i+1} + 2.
\]
$\Box$

\begin{definition}\label{D:1}
Call a partition   {\em type 1}
if successive parts differ by at least 2.
\end{definition}

\begin{corollary}
\label{type1}
For any partition $\pi$,
$\alpha(\pi)$ is a type 1 partition.
\end{corollary}

\noindent
{\bf Proof.}
By Lemma \ref{diff}, $\alpha_i - \alpha_{i+1} \geq 2 + |r_i - r_{i+1}| \geq 2$.
$\Box$

\begin{remark} The number of
partitions of $n$ with all successive ranks in $[0,1]$
is equal to the number of
type 1 partitions of $n$.
The map $\pi \rightarrow \alpha(\pi)$ is a bijection.
\end{remark}

\begin{remark} The number of
partitions of $n$ with all successive ranks in $[1,2]$
is equal to the number of
type 1 partitions of $n$ in which every part is larger than 1.
The map $\pi \rightarrow \alpha(\pi)$ is again a bijection.
\end{remark}

We can now state Andrews'
generalization of the Rogers-Ramanujan identities.  The theorem below
was established by Andrews for odd moduli $M$ \cite{A3} and was generalized
to even moduli by Bressoud \cite{B3}.

\begin{theorem}
\label{AB}
{\rm [Andrews-Bressoud]}
For integers $M$, $r$, satisfying $0 < r \leq M/2$,
let $A_n(M,r)$ be the set
of partitions of $n$ with all successive  ranks  in the interval
$[-r+2, M-r-2]$.
Then $|A_n(M,r)|$ is equal to the number
of partitions of $n$ with no part congruent to 0, $r$, or $-r$
modulo $M$.
\end{theorem}
\noindent
The smallest $M$ in the above theorem is $4$, as
when $M=3$ and  $r=1$, the interval $[-r+2, M-r-2]$ and the set $A_n(M,r)$
are both empty.
To see that Theorem \ref{AB} generalizes the Rogers-Ramanujan identities,
note that
when $r=2$  and $M=5$, the theorem says that the number
of  partitions of $n$ with all ranks in $[0,1]$ is equal to the number of
partitions of $n$ using
no part  congruent to 0, 2 or 3 modulo 5.
By Remark 1 and Theorem \ref{gordon}, this is the first Rogers-Ramanujan
identity
(\ref{RR1}).
When $r=1, M=5$, the Andrews-Bressoud theorem says that the number
of  partitions of $n$ with all ranks in $[1,2]$ is equal to the number of
partitions of $n$ using
no part congruent to 0, 1 or 4 modulo 5.
By Remark 2 and Theorem \ref{gordon}, this is the second Rogers-Ramanujan
identity (\ref{RR2}).

In fact, what Andrews and Bressoud prove is:
\begin{equation}
\label{qinfform}
\sum_{n=0}^{\infty}|A_n(M,r)|q^n \ \  = \ \
\frac1{(q)_{\infty}}\sum_{j=-\infty}^{\infty}
(-1)^jq^{j[(M)j+M-2r]/2} \ \ =
\prod_{n\not \equiv 0, \pm r (\bmod
M)}(1-q^n)^{-1}
\end{equation}
where the first equality (the hard part) follows by a sieve
argument and the second equality by application of the
Jacobi Triple Product identity.

We now show that the partitions defined by  Andrews' rank conditions
are equivalent to certain classes of multi-color partitions.

\begin{theorem}
\label{main}
For integers $r$, $M$, and $k$ satisfying
$0 < r \leq M/2$ and
$k = \lfloor M/2 \rfloor$,
let $C_n(M,r)$ be the set of
$k-1$-color partitions
$(\alpha,c_{\alpha})$ of $n$
satisfying the
following three conditions.
Let
$\alpha = (\alpha_1, \alpha_2, \ldots, \alpha_l)$.

{\rm (i) (Initial Conditions)}
For $ 1 \leq i \leq l$,
\[
\alpha_i   > \left\{ \begin{array}{ll}
	|2c_{\alpha}(i)-r+1| & \mbox{{\rm if} $\alpha_i \equiv r (\bmod{2})$}\\
	|2c_{\alpha}(i)-r| & \mbox{{\rm otherwise}}
	\end{array}
	\right. .
\]

{\rm (ii) (Color Difference Conditions)}
For $1 \leq i < l$,
\[
\alpha_i - \alpha_{i+1} \geq \left\{ \begin{array}{ll}
	2 + |2(c_{\alpha}(i) - c_{\alpha}(i+1))| & \mbox{{\rm if} $\alpha_i
\equiv
		\alpha_{i+1} (\bmod{2})$}\\
	2 + |2(c_{\alpha}(i) - c_{\alpha}(i+1))-1| & \mbox{{\rm if} $\alpha_i
\not\equiv
		\alpha_{i+1} \equiv r (\bmod{2})$}\\
	2 + |2(c_{\alpha}(i) - c_{\alpha}(i+1))+1| & \mbox{{\rm if}
$\alpha_{i+1} \not\equiv
		\alpha_{i} \equiv r (\bmod{2})$}
	\end{array}
        \right. .
\]

{\rm (iii) (Parity Condition on Last Color when $M$ is Even)}
For $1 \leq i \leq l$, if $M$ is even and $c_{\alpha}(i) = k-1$,
the last color, then
\[
\alpha_i \not\equiv r \pmod 2.
\]
Then the number of partitions in the two sets $A_n(M, r)$
and  $C_n(M, r)$ are the same.

\end{theorem}

\noindent
{\bf Proof.}
This follows from (\ref{qinfform}) once we establish a bijection
between $A_n(M,r)$ and $C_n(M,r)$.
Let $\pi$ be a partition of $n$ with ranks
$r_1, r_2, \ldots r_{d(\pi)}$, and assume that all ranks lie  in the interval
$[-r+2,M-r-2]$.
For $1 \leq i \leq d$, let
$\alpha_i = \pi_i + \pi_{i}' -2i+1$, as before,
and let $\alpha(\pi)$ be the partition of $n$ defined by
$\alpha_1,
\alpha_2, \ldots,
\alpha_d$.

Color the parts of $\alpha$ as follows:  for $1 \leq i \leq d$,
\begin{equation}
c_{\alpha}(i) =  \left\{ \begin{array}{ll}
	(r_i + r - 1)/2 &  \mbox{{\rm if} $\alpha_i \equiv r (\bmod{2})$}\\
	(r_i + r )/2 &  \mbox{{\rm otherwise}}.
	\end{array}
	\right.
\label{assigncolor}
\end{equation}
We show that the $k-1$-color partition $(\alpha, c_{\alpha})$
satisfies conditions (i)-(iii) of the theorem.

\noindent
Condition(i).
By Lemma \ref{partsize},
$\alpha_i > |r_i|$
and from (\ref{assigncolor}),
\begin{equation}
r_i =  \left\{ \begin{array}{ll}
	2c_{\alpha}(i)-r+1 &  \mbox{{\rm if} $\alpha_i \equiv r (\bmod{2})$}\\
	2c_{\alpha}(i)-r &  \mbox{{\rm otherwise}}.
	\end{array}
        \right.
\label{rankfromcolor}
\end{equation}

\noindent
Condition(ii).
By Lemma \ref{diff}, $\alpha_i - \alpha_{i+1} \geq 2 + |r_i - r_{i+1}|$
and the conditions follow from (\ref{rankfromcolor}).

\noindent
Condition(iii).
If $M=2k$ and $c_{\alpha}(i)=k-1$, then by definition of
$c_{\alpha}$,
\[
r_i \in \{-r+2k -2, \ -r+2k-1\} = \{M-r-2,M-r-1\}.
\]
But
$r_i \not =  M-r-1$ since $r_i \in [-r+2, M-r-2]$.
It follows then that $r_i = M-r-2$ and therefore that
\begin{eqnarray*}
\alpha_i  &  =  &  \pi_i + \pi'_i -2i+1=r_i + 2\pi'_i - 2i + 1\\
          &  =  &  M-r -1 +2(\pi'_i - i)\equiv r +1 \pmod 2.
\end{eqnarray*}

To show this is a bijection, given $r, M$, and $k$ satisfying the conditions
of the theorem, let
$(\alpha,c_{\alpha})$ be a $k-1$ color partition of $n$, with
$\alpha = (\alpha_1, \ldots, \alpha_l)$,
satisfying (i), (ii), and (iii) of the theorem.
We define an inverse map which sends $\alpha$ to a partition $\pi$ of $n$
with $d$ angles, where the $i$th angle of $\pi$ has width $x_i$ defined by
\begin{equation}\label{E:invmap}
x_i =  \left\{  \begin{array}{ll}
	(-r+2c_{\alpha}(i)+\alpha_i + 2)/2 & \mbox{{\rm if} $\alpha_i \equiv r
(\bmod{2})$}\\
	(-r+2c_{\alpha}(i)+\alpha_i + 1)/2 & \mbox{{\rm otherwise}}
	\end{array}
	\right.
\label{inverse}
\end{equation}
and height $y_i=\alpha_i-x_i+1$.
We must verify that $\pi$ is a partition, i.e. that
$x_1 > x_2 > \ldots > x_d \geq 1$ and
$y_1 > y_2 > \ldots > y_d \geq 1$ (clearly it has weight $n$),
and that the $i$th rank $r_i=x_i-y_i$ lies in the interval
$[-r+2,M-r-2]$ and furthermore that it satisfies (\ref{assigncolor}).

We first verify that for $1 \leq i \leq l$, $x_i \geq 1$ and $y_i \geq 1$.
If $\alpha_i \equiv r (\bmod{2})$, then using (\ref{inverse}) and condition
(i),
\[
x_i = (-r+2c_{\alpha}(i)+\alpha_i+2)/2 >
(-r+2c_{\alpha}(i)+|-2c_{\alpha}(i)+r-1|+2)/2\geq 1/2
\]
and
\[
y_i = (\alpha_i+r-2c_{\alpha}(i))/2>(|2c_{\alpha}(i)-r+1|+r-2c_{\alpha}(i))/2
\geq 1/2.
\]
Similarly,
if $\alpha_i \not\equiv r (\bmod{2})$,
\[
x_i = (-r+2c_{\alpha}(i)+\alpha_i+1)/2  >
(-r+2c_{\alpha}(i)+|-2c_{\alpha}(i)+r|+1)/2\geq  1/2
\]
and
\[
y_i = (\alpha_i+r-2c_{\alpha}(i)+1)/2 >
(|2c_{\alpha}(i)-r|+r-2c_{\alpha}(i)+1)/2
\geq  1/2.
\]
Checking the rank of the $i$th angle
\[
x_i-y_i = 2x_i-\alpha_i-1=\left\{ \begin{array}{ll}
	-r+2c_{\alpha}(i)+1 & \mbox{{\rm if} $\alpha_i \equiv r (\bmod{2})$}\\
	-r+2c_{\alpha}(i) & \mbox{{\rm otherwise}}
	\end{array},
	\right.
\]
which satisfies (\ref{assigncolor}).

Now we verify that for $1 \leq i < l$, $x_i > x_{i+1}$ and
$y_i > y_{i+1}$.
We check each of three cases using (\ref{inverse}) and condition
(ii) of the theorem.
If $\alpha_i \equiv \alpha_{i+1} (\bmod{2})$, then
\[
x_i - x_{i+1} =(\alpha_i-\alpha_{i+1} + 2(c_{\alpha}(i)-c_{\alpha}(i+1)))/2
\geq 1,
\]
\[
y_i - y_{i+1} =(\alpha_i-\alpha_{i+1} - 2(c_{\alpha}(i)-c_{\alpha}(i+1)))/2
\geq 1.
\]
If $\alpha_i \not\equiv \alpha_{i+1} \equiv r (\bmod{2})$, then
\[
x_i - x_{i+1} =(\alpha_i-\alpha_{i+1} + 2(c_{\alpha}(i)-c_{\alpha}(i+1))-1)/2
\geq 1,
\]
\[
y_i - y_{i+1} =(\alpha_i-\alpha_{i+1} - 2(c_{\alpha}(i)-c_{\alpha}(i+1))+1)/2
\geq 1.
\]
Finally, if $\alpha_{i+1} \not\equiv \alpha_{i}  \equiv r (\bmod{2})$, then
\[
x_i - x_{i+1} =(\alpha_i-\alpha_{i+1} + 2(c_{\alpha}(i)-c_{\alpha}(i+1))+1)/2
\geq 1,
\]
\[
y_i - y_{i+1} =(\alpha_i-\alpha_{i+1} - 2(c_{\alpha}(i)-c_{\alpha}(i+1))-1)/2
\geq 1.
\]

$\Box$

The following result is an immediate consequence of Theorem \ref{main} and
(\ref{qinfform}).

\begin{corollary}

\[
\sum_{n=0}^{\infty}|C_n(M,r)|q^n \ \  = \ \
\frac1{(q)_{\infty}}\sum_{j=-\infty}^{\infty}
(-1)^jq^{j[(M)j+M-2r]/2} \ \ =
\prod_{n\not \equiv 0, \pm r (\bmod
M)}(1-q^n)^{-1}.
\]

\end{corollary}

$\Box$

See Figure \ref{example1} for an example of the bijection when $M=7,r=1$, and
$n=10$.
Figure \ref{example2} shows the bijection when $M=8,r=3$, and $n=10$.

\section{Remarks on Alternative Colorings}

An alternative coloring is to color the part $\alpha_i$, derived from the
$i$th angle with rank $r_i$, by:
\begin{equation}
c'_{\alpha}(i) = \left\{ \begin{array}{ll}
	r_i-k+r & \mbox{{\rm if} $r_i > k-r$}\\
	-r_i+k-r & \mbox{{\rm otherwise} }
	\end{array}
	\right.
\label{color2}
\end{equation}
and proceed to formulate the conditions required to make this
mapping a bijection between $(k-1)$-color partitions of $n$
which satisfy the conditions and the set $A_{n}(M,r)$
of partitions of $n$ with all ranks in $[-r+2,M-r-1]$.
This is the coloring which would follow by applying coloring
ideas in Agarwal and Andrews (for other families)
\cite{AgAn} to the partition family
$A_n(M,r)$.
This same multi-coloring, $c'_{\alpha}$, arises by applying the
coloring scheme in Agarwal and Bressoud \cite{AgBr}
(where it was used on a different
family) to the lattice path
interpretation of
$A_n(M,r)$
described in Bressoud \cite{B4},
namely, for each peak of height $y$ at location $x$ in the lattice
path, create a part of size $x$ with color $y$.
This can be shown to be equivalent to the
coloring (\ref{color2})
using the bijection between the lattice path and rank interpretations
described by Bressoud in \cite{B4}.
An advantage of our coloring is that the
color depends only on $r$ and $r_i$ and {\em not on $k$}.
This may make it easier
to establish a direct connection between
the multi-sum and our multi-color interpretation.
The hope is that for our particular goals,
this interpretation will be more fruitful in
the Lie algebra setting.

\section{Finitized Rogers-Ramanujan Identities}

In view of (\ref{verygenRR}) and the second equality
in (\ref{qinfform}),
the generalized Rogers-Ramanujan identities can be written as
\begin{equation}
\label{physicsform}
\frac1{(q)_{\infty}}\sum_{j=-\infty}^{\infty}
(-1)^jq^{j[(M)j+M-2r]/2} \ \ =
\sum_{n_1 \geq  \ldots \geq n_{k-1}\geq 0}
\frac{q^{n_1^2+n_2^2+ \ldots + n_{k-1}^2 + n_r + n_{r+1} + \ldots + n_{k-1}}}
{(q)_{n_1-n_2} \ldots (q)_{n_{k-2}-n_{k-1}}(q^{2-s}; q^{2-s})_{n_{k-1}}},
\end{equation}
where $s=0$ if $M$ is even and $s=1$ if $M$ is odd.
In fact, this is the form of interest in recent work
relating Rogers-Ramanujan identities to applications
in statistical mechanics and conformal field theory.
In these applications, the left-hand side is the
{\em bosonic form} and the right-hand side is the
{\em fermionic form}.
In the general case, the bosonic form is associated with a
character of the minimal model of a Virasoro algebra.
In \cite{LW1a} and later in \cite{LP} a
generalized fermionic Pauli exclusion principle
was discovered and discussed
in connection with higher level representations of the affine
Lie algebra $\widehat{\mathfrak{sl}}_2$. Using similar ideas the sum
side of the generalized Rogers-Ramanujan identities has recently been studied
in several papers, and it has been shown that under certain restrictions,
every bosonic form has a fermionic form and, even more,
there is a corresponding finitization \cite{BMS}.

In \cite{FQ}, motivated by these connections to physics,
Foda and Quano derived
finite approximations of the
identities (\ref{physicsform}).
We will
show that the multi-color partition interpretation also
provides an interpretation of the left-hand side for the
finitization in this case.

We first make an observation.
For $0 < r \leq M/2$ and $u,v \geq 0$,
let $F_n(M, r, u, v)$ be the set of partitions of $n$
with all ranks in the interval $[-r+2, M-r-2]$ and
whose Ferrer diagrams are
contained in a $v \times u$ rectangle.

\begin{lemma}
\label{box}
For $0 < r \leq M/2$ and $u,v \geq 0$,
the set $F_n(M, r, u, v)$ is in one-to-one correspondence
with the set of
$k-1$-color partitions
$(\alpha,c_{\alpha})$ of $n$ in $C_n(M,r)$
with the following additional constraint
\[
(iv) \ \ (u+v-1)-\alpha_1 \geq
	|2c_{\alpha}(1)-r+(v-u)+\frac 12|-\frac 12.
\]
\end{lemma}
\noindent {\bf Proof.} Under our map
in the proof of Theorem \ref{main},  the partitions in
$F_n(M, r, u, v)$ are mapped into $k-1$ color
partitions $(\alpha_i)$ satisfying (i), (ii), and (iii) in Theorem \ref{main}.
Since this map was shown to be a bijection, it suffices to show that
its
image is characterized by the extra condition (iv).
So, for a $k-1$ color partition $\alpha$ satisfying (i)-(iii), we
derive necessary and sufficient
conditions on $\alpha$ to guarantee that under the inverse mapping
its first angle has width $x_1\leq u$
and  height $y_1\leq v$.
The inverse map (\ref{E:invmap}) has
\[
x_1 = \lfloor (-r+2c_{\alpha}(1)+\alpha_1+ 2)/2 \rfloor
\]
so $x_1 \leq u$ is equivalent to
\[
\frac{-r+2c_{\alpha}(1)+\alpha_1+ 1}2 \leq u
\]
which in turn is equivalent to
\begin{equation}
\label{ubound}
2c_{\alpha}(1)-r + (v-u) \leq (u+v-1) - \alpha_1.
\end{equation}
On the other hand, the height $y_1=\alpha_1-x_1+1$
must satisfy $y_1 \leq v$ which is equivalent to
$x_1 \geq -v+\alpha_1+1$, so we require
\[
\frac{-r+2c_{\alpha}(1)+\alpha_1+2}2\geq -v+\alpha_1+1
\]
which
is equivalent to
\begin{equation}
\label{vbound}
-2c_{\alpha}(1)+r - (v-u) -1 \leq (u+v-1) - \alpha_1.
\end{equation}
Combining (\ref{ubound}) and (\ref{vbound}) gives the
lemma.
$\Box$

In \cite{FQ}, Foda and Quano derive the following
finitized Rogers-Ramanujan identity:
\begin{eqnarray}
\label{finitizeodd}
\sum_{j=-\infty}^{\infty}&&
(-1)^jq^{j[(2k+1)j+2k+1-2r]/2}\left[\begin{array}{c} N \\
\lfloor\frac{N-k+r-(2k+1)j}2\rfloor\end{array}\right]_q \\
&=&\sum q^{n_1^2+\cdots+n_{k-1}^2+n_{r}+\cdots+n_{k-1}}
\prod_{j=1}^{k-1} \left[\begin{array}{c}
N-2(n_1+\cdots+n_{j-1})-n_j-n_{j+1}-a_{rj}^{(k)}\\
n_j-n_{j+1}\end{array}\right]_q, \nonumber
\end{eqnarray}
where the sum runs over $n_1\geq \cdots\geq n_{j}=0$ such that $
2(n_1+\cdots+n_{k-1})\leq N-k+r$ and
$a_{ij}^{(k)}$ is the $ij$-entry of the $k\times (k-1)$ matrix
\[
A^{(k)}=\left(\begin{array}{cccc} 1 & 2 & \cdots & k-1\\
0 & 1 & \cdots & k-2\\ \vdots & \vdots &\ddots & \vdots \\
0 & 0 &\cdots & 1 \\ 0 & 0 &\cdots & 0\end{array}\right).
\]
When $N\rightarrow \infty$, the identity reduces to
(\ref{physicsform}) in the case of odd $M$ and $s=1$.

Foda and Quano prove that the left-hand side of (\ref{finitizeodd}) is the
generating function for the set of partitions
$F_n(M,r,u,v)$ in the special case that
$M=2k+1$,
$u = \lfloor    (N+k-r+1)/2 \rfloor$, and
$v = \lfloor    (N-k+r)/2 \rfloor$.
Combining this with Lemma (\ref{box}) we get the following.
\begin{corollary}
For $M=2k+1$ and $r \leq k$,
the left-hand side of (\ref{finitizeodd}) is the generating function
for the set of
$k-1$-color partitions
$(\alpha,c_{\alpha})$ of $n$ in $C_n(M,r)$
with the following additional constraint
\[
(iv) \ \ (N-1)-\alpha_1 \geq
\left\{ \begin{array}{ll}
	|2c_{\alpha}(1)-k+\frac 12|-\frac 12, & \mbox{{\rm if} $N+k \equiv
r (\bmod {2})$}\\
	|2c_{\alpha}(1)-k-\frac 12|-\frac 12, & \mbox{{\rm otherwise}}

	\end{array}
	\right.
\]
In particular, $\alpha_1\leq N-1$.
\end{corollary}

\pagebreak
In \cite{FQ}, Foda and Quano  derive the following additional
finitized Rogers-Ramanujan identity:
\begin{eqnarray}
\label{finitizeeven}
\sum_{j=-\infty}^{\infty}
(-1)^jq^{j(kj+k-r)}\left[\begin{array}{c} 2N+k-r \\
N-kj \end{array}\right]_q =
\sum q^{n_1^2+\cdots+n_{k-1}^2+n_{r}+\cdots+n_{k-1}} \times\\
\prod_{j=1}^{k-2} \left[\begin{array}{c}
2N-2(n_1+\cdots+n_{j-1})-n_j-n_{j+1}-b_{rj}^{(k)}\\
n_j-n_{j+1}\end{array}\right]_q
\left[\begin{array}{c}
N-(n_1+\cdots+n_{k-2})\\
n_{k-1}\end{array}\right]_{q^2},\nonumber
\end{eqnarray}
where the sum runs over $n_1\geq \cdots\geq n_{j}=0$ such that $
n_1+\cdots+n_{k-1}\leq N$ and
$b_{ij}^{(k)}$ is the $ij$-entry of the $k\times (k-2)$ matrix
\[
B^{(k)}=\left(\begin{array}{ccccc}
k-2 & k-3 & k-4 &\cdots & 1\\
k-2 & k-3 & k-4 &\cdots & 1\\
k-3 & k-4 & k-5 &\cdots & 1\\
\vdots &\vdots  &\vdots  &\cdots &\vdots \\
3 & 2 & 1 & \cdots & 1\\
2 & 1 & 1 & \cdots  &1\\
1 & 1 & 1 & \cdots  &1\\
0 & 0 & 0 & \cdots  &0\\
\end{array}\right).
\]
When $N\rightarrow \infty$, the identity reduces to
(\ref{physicsform}) in the case of even $M$ and $s=0$.

Foda and Quano prove that the left-hand side of (\ref{finitizeeven}) is the
generating function for the set of partitions
$F_n(M,r,u,v)$ in the special case that
$M=2k$,
$u = N+k-r$, and
$v = N$.
Combining this with Lemma \ref{box} we get:
\begin{corollary}
For $M=2k$, the left-hand side of (\ref{finitizeeven})
is the generating function
for the set of
$k-1$-color partitions
$(\alpha,c_{\alpha})$ of $n$ in $C_n(M,r)$
and the following additional constraint
\[
(iv) \ \ (2N+k-r-1)-\alpha_1 \geq
	|2c_{\alpha}(1)-k+\frac 12|-\frac 12.
\]
\end{corollary}

\begin{remark}
Clearly, our map in Theorem \ref{main}, when restricted to partitions in
$F_n(M, r, u,v)$ with Durfee square of size $d$, gives a bijection
with those $k-1$ color partitions of $n$ in $C_n(M,r)$
which have exactly $d$ parts and satisfy (iv) of
of Lemma \ref{box}.
\end{remark}

\vspace{.2in}
\noindent
{\bf Acknowledgment}
We would like to thank James Lepowsky for helpful comments on an
earlier version of this paper and Christian Krattenthaler
for bringing the references \cite{AgAn} and \cite{AgBr} to our
attention.

\begin{figure}
{\small
\begin{tabbing}
xxxxxxxxxxxxxxx\=xxxxxxxxxxxxxxxxxx\=xxxxxxxxxxxxxxxxxxxxxx\= \kill
\>Partition  $\pi$ \>   rank vector\>	the 2-color partition\\
\>of 10 with all\>      of $\pi$\>   	$(\alpha(\pi),c_{\alpha(\pi)})$ of 10\\
 \>ranks in $[1,4]$\\
\>\\
\>$(7,1,1,1)$\>	$[3]$\>	$(10_2)$\\
\>$(6,4)$	\>	$[4,2]$\>	$(7_2,3_1)$\\
\>$(6,3,1)$\>	$[3,1]$ \>$(8_2,2_1)$\\
\>$(6,1,1,1)$\>	$[1]$\>	$(10_1)$\\
\>$(5,5)$	\>	$[3,3]$\>	$(6_2,4_2)$\\
\>$(5,4,1)$\>	$[2,2]$\>	$(7_1,3_1)$\\
\>$(5,3,1,1)$\>	$[1,1]$\>	$(8_1,2_1)$\\
\>$(4,4,2)$\>	$[1,1]$\>	$(6_1,4_1)$\\
\end{tabbing}
}
\caption{Example of the bijection of Theorem \ref{main} when $M=7$, $r=1$, and
$n=10$.}

\label{example1}
\end{figure}

\begin{figure}
{\small
\begin{tabbing}
xxxxxxxxxxxxxxx\=xxxxxxxxxxxxxxxxxx\=xxxxxxxxxxxxxxxxxxxxxx\= \kill
\>Partition  $\pi$ \>   rank vector\>	the 3-color partition\\
\>of 10 with all\>      of $\pi$\>   	$(\alpha(\pi),c_{\alpha(\pi)})$ of 10\\
 \>ranks in $[-1,3]$\\
\>\\
\>$(6,2,1,1)$ \>	$[2,0]$ \>		$(9_2,1_1)$\\
\>$(6,1,1,1,1)$ \>	$[1]$ \>		$(10_2)$\\
\>$(5,4,1)$ \>	$[2,2]$ \>		$(7_2,3_2)$\\
\>$(5,3,2)$ \>	$[2,0]$ \> 	$(7_2,3_1)$\\
\>$(5,3,1,1)$ \>	$[1,1]$ \>		$(8_2,2_2)$\\
\>$(5,2,2,1)$ \>	$[1,-1]$ \> 	$(8_2,2_1)$\\
\>$(5,2,1,1,1)$ \>	$[0,0]$	 \> 	$(9_1,1_1)$\\
\>$(5,1,1,1,1,1)$ \>	$[-1]$  \>  	$(10_1)$\\
\>$(4,4,2)$ \>	$[1,1]$	  \>	$(6_2,4_2)$\\
\>$(4,4,1,1)$ \>	$[0,2]$	  \>	$(7_1,3_2)$\\
\>$(4,3,3)$ \>	$[1,0,0]$   \>	$(6_2,3_1,1_1)$\\
\>$(4,3,2,1)$ \>	$[0,0]$	 \>	$(7_1,3_1)$\\
\>$(4,3,1,1,1)$ \>	$[-1,1]$ \>	$(8_1,2_2)$\\
\>$(4,2,2,1,1)$ \>	$[-1,-1]$ \>	$(8_1,2_1)$\\
\>$(3,3,3,1)$ \>	$[-1,0,0]$ \>	$(6_1,3_1,1_1)$\\
\>$(3,3,2,2)$ \>	$[-1 -1]$ \>	$(6_1,4_1)$\\
\>$(7,1,1,1)$ \>	$[3]$	 \>	$(10_3)$\\
\>$(6,3,1)$ \>	$[3,1]$	 \>	$(8_3,2_2)$\\
\>$(6,2,2)$ \>	$[3,-1]$  \>	$(8_3,2_1)$\\
\>$(5,5)$	 \>	$[3,3]$   \>	$(6_3,4_3)$
\end{tabbing}
\caption{Example of the bijection of Theorem \ref{main} when $M=8$, $r=3$, and
$n=10$.}
\label{example2}
}
\end{figure}

\pagebreak

\end{document}